\documentclass[12pt,a4paper,oneside,reqno,notitlepage]{amsart}
\usepackage{amssymb}
\textwidth
150mm

\newtheorem{theorem}{Theorem}
\theoremstyle{plain}

\newtheorem{definition}{Definition}

\newtheorem{lemma}{Lemma}

\newtheorem{proposition}{Proposition}

\numberwithin{equation}{section}

\begin{document}

\baselineskip 8mm
\parindent 9mm

\title[]
{No local $L^{1}$ solutions for semilinear fractional heat equations}

\author{ Kexue Li}

\address{School of Mathematics and Statistics,
Xi'an Jiaotong University,
 Xi'an
710049, China; Department of Mathematics, University of Washington, Seattle, WA, 98195, USA}
\email{kexueli@gmail.com}

\thanks{{\it 2010 Mathematics Subjects Classification}: 35K55}

\begin{abstract}
We study the Cauchy problem for the semilinear fractional heat equation $u_{t}=\triangle^{\alpha/2}u+f(u)$ with non-negative initial value $u_{0}\in L^{q}(\mathbb{R}^{n})$ and locally Lipschitz, non-negative source term $f$. For $f$ satisfying the Osgood-type condition $\int_{1}^{\infty}\frac{ds}{f(s)}=\infty$, we show that there exist initial conditions such that the equation has no local solution in $L^{1}_{loc}(\mathbb{R}^{n})$.

\end{abstract}
\maketitle

\section{\textbf{Introduction}}
In this paper, we consider the Cauchy problem for the semilinear fractional heat equation
\begin{align}\label{f}
u_{t}=\triangle^{\alpha/2}u+f(u), \ u(0)=u_{0}\in L^{q}(\mathbb{R}^{n}), \ q\geq1,
\end{align}
where $\triangle^{\alpha/2}:=-(-\Delta)^{\alpha/2}$ denotes the fractional Laplacian defined by the Fourier transform
\begin{align*}
(\mathcal{F}(-\triangle)^{\alpha/2}u)(\xi)=|\xi|^{\alpha}\mathcal{F}(u)(\xi),
\end{align*}
where $0<\alpha\leq2$, $\mathcal{F}$ denotes the Fourier transform. \\
We assume the following conditions hold for the initial value $u_{0}$ and source term $f$:  \\
(A1) $u_{0}\geq 0$ and  $u_{0}\in L^{q}(\mathbb{R}^{n})$;\\
(A2) $f: [0,\infty)\rightarrow [0,\infty)$ is locally Lipschitz continuous, non-decreasing, $f(0)=0$, $f>0$ on $(0,\infty)$;\\
(A3) $f$ satisfies the Osgood-type condition
\begin{align}\label{Osgood}
\int_{1}^{\infty}\frac{ds}{f(s)}=\infty.
\end{align}

Recall ordinary differential equations of the form
\begin{align}\label{ordinary}
u_{t}=f(u),
\end{align}
where $f$ is positive, continuous and satisfies (\ref{Osgood}).
It is well known that this Osgood-type condition is necessary and sufficient for the global existence of solutions of (\ref{ordinary}) with initial data $u_{0}\geq 0$.

For $\alpha=2$, (\ref{f})  reduces to
\begin{align}\label{Laplacian}
u_{t}=\triangle u+f(u), \ u(0)=u_{0}\in \mathbb{R}^{n}.
\end{align}
When $f$ is locally Lipschitz continuous, by Theorem 1.4 in \cite{Pazy}, there exists a $t_{max}\leq \infty$ such that a unique mild solution of (\ref{Laplacian}) exists on $[0,t_{max})$, if $t_{max}<\infty$, then $\lim_{t\uparrow t_{max}}\|u(t)\|=\infty$. For non-negative $u_{0}\in L^{\infty}(\mathbb{R}^{n})$ and $f$ satisfying  the Osgood-type condition (\ref{Osgood}), we can get the global existence of (\ref{Laplacian}) (See, for example, Remark 5.1 in \cite{HB}). However, when the initial data $u_{0}$ is singular or unbounded, the question whether (\ref{Laplacian}) has global existence hasn't been settled until the appearance of \cite{LRS}. Existence results for (\ref{Laplacian}) is considered in \cite{LRS} under the conditions that the source term $f$ satisfies conditions (A2) and (A3). The authors showed that there are initial conditions satisfying (A1) for which there is no local integral solution (See the following Definition \ref{local integral solution}) of (\ref{Laplacian}) that remains in $L_{loc}^{1}(\mathbb{R}^{n})$.
\begin{definition}\label{local integral solution}
(See \cite{QS}, p.78) Given $f$ non-negative and $u_{0}\geq 0$, we say that $u$ is a local integral solution of (\ref{Laplacian}) on $[0,T)$ if $u: \mathbb{R}^{n}\times [0,T)\rightarrow [0,\infty]$ is measurable, finite almost everywhere, and
\begin{align*}
u(t)=S(t)u_{0}+\int_{0}^{t}S(t-s)f(u(s))ds
\end{align*}
holds almost everywhere on $\mathbb{R}^{n}\times [0,T)$, where
\begin{align*}
(S(t)u_{0})(x)=(4\pi t)^{-n/2}\int_{\mathbb{R}^{n}}e^{-|x-y|^{2}/4t}u_{0}(y)dy
\end{align*}
is a classical solution of the linear heat equation
\begin{align*}
v_{t}=\Delta v,\ v(0)=u_{0}.
\end{align*}
\end{definition}

From Definition \ref{local integral solution}, we see that if $u$ is a classical solution (or mild solution), then it is an integral solution. Non-existence integral solutions implies the non-existence of classical solutions (or mild solutions). The non-existence of a local solution is `instantaneous blow-up' in some sense (see, e.g., \cite{PV}).

In fact, the Osgood-type condition (\ref{Osgood}) is not necessary for global existence of solutions for (\ref{Laplacian}). Fujita \cite{Fu,Fujita} studied the initial value problem
\begin{align}\label{p}
u_{t}=\triangle u+u^{p}, \ u(0)=u_{0}\in \mathbb{R}^{n},
\end{align}
where $u_{0}\geq 0$, $p>1$. He proved that if $p>2/n+1$, then global solutions exists for small initial data $u_{0}$. Weissler \cite{Weissler} considered the problem (\ref{p}), in the case $p>2/n+1$, $\|u_{0}\|_{L^{n(p-1)/2}}$ is sufficiently small, the global solution is obtained. It is clear that the source term $f(u)=u^{p}(p>1)$ doesn't satisfy the Osgood-type condition (\ref{Osgood}).

In recent years, a great deal of attention has been paid to fractional differential equations with fractional Laplacian due to their many applications in mathematics, physics, biology, see for instance \cite{AT,MDH,SJ,CGL,GH} and references therein. It is nature to consider the question whether the Osgood-type condition (\ref{Osgood}) guarantees the global existence of solutions for semilinear fractional heat equations (\ref{f}) with singular (or unbounded) initial data. We answer this question in the negative.

The paper is organized as follows. In Section 2, we present some estimates for the linear fractional heat equation. In Section 3, we construct a function which satisfies the Osgood-type condition. In Section 4, we prove the non-existence results for  (\ref{f}) with unbounded initial data.

\section{\textbf{Estimates for some solutions of the fractional heat equation}}
For any $r>0$, let $B_{r}(x)$ be the Euclidean ball in $\mathbb{R}^{n}$ of radius $R$ centred at $x$, $S^{n-1}$ the unit sphere and $\omega_{n}$ the volume of the unit ball $B_{1}(0)$.
Let $p(t,x,y)$ be the heat kernel of $\Delta ^{\alpha/2}$ on $\mathbb{R}^{n}$.
Throughout this paper, we use $c_{0}$, $c_{1}$, $c_{2}$, $\ldots$ to denote generic constants, which may change from line to line.
For two nonnegative functions $f_{1}$ and $f_{2}$, the notion $f_{1}\asymp f_{2}$ means that $c_{1}f_{2}(x)\leq f_{1}(x)\leq c_{2}f_{2}(x)$, where $c_{1},c_{2}$ are positive constants.
It is well known that (see, e.g., \cite{TG,KT,CKS})
\begin{align*}
p(t,x,y)\asymp \big(t^{-n/\alpha}\wedge \frac{t}{|x-y|^{n+\alpha}}\big),
\end{align*}
that is, there exist constants $c_{1}$, $c_{2}$ such that for $t>0$, $x,y\in \mathbb{R}^{n}$,
\begin{align}\label{c1}
c_{1}\big(t^{-n/\alpha}\wedge \frac{t}{|x-y|^{n+\alpha}}\big)\leq p(t,x,y)\leq c_{2}\big(t^{-n/\alpha}\wedge \frac{t}{|x-y|^{n+\alpha}}\big),
\end{align}
where $c_{1}$ and $c_{2}$ are positive constants depending on $\alpha$. \\
From (\ref{c1}), it follows
\begin{align}\label{alpha}
\frac{c_{3}t}{(t^{1/\alpha}+|y-x|)^{n+\alpha}}\leq p(t,x,y)\leq \frac{c_{4}t}{(t^{1/\alpha}+|y-x|)^{n+\alpha}},
\end{align}
where $t>0$, $x,y\in \mathbb{R}^{n}$, $c_{3}$ and $c_{4}$ are positive constants depending on $\alpha$.
\begin{definition}
$u$ is said to be a local integral solution of (\ref{f}) on $[0,T)$, if $u: \mathbb{R}^{n}\times [0,T)\rightarrow [0,\infty)$ is measurable, finite almost everywhere and
\begin{align}\label{solution}
u(t)=S_{\alpha}(t)u_{0}+\int_{0}^{t}S_{\alpha}(t-s)f(u(s))ds
\end{align}
holds almost everywhere in $\mathbb{R}^{n}\times [0,T)$, where $(S_{\alpha}(t)u_{0})(x)$ is a classical solution of the linear fractional heat equation
\begin{align}\label{linear}
w_{t}=\triangle^{\alpha/2}w, \ w(0)=u_{0},
\end{align}
and
\begin{align}\label{fractional}
(S_{\alpha}(t)u_{0})(x)=\int_{\mathbb{R}^{n}}p(t,x,y)u_{0}(y)dy,
\end{align}
where $p(t,x,y)$ is the heat kernel of $\Delta^{\alpha/2}$ on ${R}^{n}$.
\end{definition}
\begin{proposition}\label{volume}
Let $\beta\in (0,n)$ and $R>1$. Assume that $u_{0}\in L^{1}(\mathbb{R}^{n})$ be the non-negative, radially symmetric function given by
\begin{equation*}
u_{0}(x)=|x|^{-\beta}\chi_{R}:=\left\{\begin{aligned}
&|x|^{-\beta}, \ |x|\leq R,\\
&0, \ \ \ \ \ \ |x|>R.
\end{aligned}\right.
\end{equation*}
Let $w(t)=S_{\alpha}(t)u_{0}$ and
\begin{align}\label{m}
M=\min\{w(\tilde{x},t): \tilde{x}\in S^{n-1}, \ 0\leq t\leq 1\}.
\end{align}
If $\gamma\in (0, 1/\alpha)$, then for any $\phi\geq c_{3}M/c_{4}$,
\begin{align}
w(x,t)\geq \phi \ for \ |x|\leq t^{\gamma} \ and \ 0<t\leq (c_{4}\phi/c_{3}M)^{-1/\beta\gamma},
\end{align}
where $c_{3}$,  $c_{4}$ are the same constants as in (\ref{alpha}).
\end{proposition}

\begin{proof}
For any $0<t\leq 1$, $x=t^{\gamma}\tilde{x}\in \partial B_{n}(t^{\gamma})$, $\tilde{x}\in S^{n-1}$,
\begin{align}\label{b}
w(x,t)&=w(t^{\gamma}\tilde{x},t)\nonumber\\
&\geq \frac{c_{3}}{t^{n/\alpha}}\int_{\mathbb{R}^{n}}\frac{u_{0}(y)}{(1+t^{-1/\alpha}|y-x|)^{n+\alpha}}dy\nonumber\\
&=\frac{c_{3}}{t^{n/\alpha}}\int_{\mathbb{R}^{n}}\frac{u_{0}(y)}{(1+t^{-1/\alpha}|y-t^{\gamma}\tilde{x}|)^{n+\alpha}}dy\nonumber\\
&=\frac{c_{3}}{t^{n/\alpha}}\int_{\mathbb{R}^{n}}\frac{u_{0}(t^{\gamma }z)t^{n \gamma}}{(1+t^{\gamma-1/\alpha}|z-\tilde{x}|)^{n+\alpha}}dz\nonumber\\
&=\frac{c_{3}}{t^{n/\alpha}}\int_{B_{t^{-\gamma}R}(0)}\frac{|t^{\gamma}z|^{-\beta}t^{n\gamma}}{(1+t^{\gamma-1/\alpha}|z-\tilde{x}|)^{n+\alpha}}dz\nonumber\\
&=\frac{t^{-\beta\gamma}c_{3}}{t^{n/\alpha}}\int_{B_{t^{-\gamma}R}(0)}\frac{|z|^{-\beta}t^{n\gamma}}{(1+t^{\gamma-1/\alpha}|z-\tilde{x}|)^{n+\alpha}}dz\nonumber\\
&= \frac{t^{-\beta\gamma}c_{3}}{t^{n(\frac{1}{\alpha}-\gamma)}}\int_{B_{t^{-\gamma}R}(0)}\frac{|z|^{-\beta}}{(1+t^{\gamma-1/\alpha}|z-\tilde{x}|)^{n+\alpha}}dz\nonumber\\
&\geq \frac{t^{-\beta\gamma}c_{3}}{t^{n(\frac{1}{\alpha}-\gamma)}}\int_{B_{R}(0)}\frac{|z|^{-\beta}}{(1+t^{\gamma-1/\alpha}|z-\tilde{x}|)^{n+\alpha}}dz\nonumber\\
&=\frac{t^{-\beta\gamma}c_{3}}{t^{n(\frac{1}{\alpha}-\gamma)}}\int_{\mathbb{R}^{n}}\frac{u_{0}(z)}{(1+t^{\gamma-1/\alpha}|z-\tilde{x}|)^{n+\alpha}}dz.
\end{align}
Since $\gamma\in (0,1/\alpha)$,  we have $0<\gamma\alpha<1$. Then from (\ref{b}), it follows that
\begin{align}\label{lower bound}
w(x,t)\geq  \frac{t^{-\beta\gamma}c_{3}}{t^{n(\frac{1}{\alpha}-\gamma)}}\int_{R^{n}}\frac{u_{0}(z)}{(1+t^{\gamma-1/\alpha}|z-\tilde{x}|)^{n+\frac{1}{\frac{1}{\alpha}-\gamma}}}dz.
\end{align}
This together with (\ref{fractional}) and (\ref{alpha}) yield
\begin{align}
w(x,t)\geq \frac{c_{3}}{c_{4}}t^{-\beta\gamma}w(\tilde{x},t^{\alpha/(1-\alpha \gamma)}).
\end{align}
For $0<t\leq 1$ and $|\tilde{x}|=1$, by (\ref{m}),
\begin{align*}
w(\tilde{x},t^{\alpha/(1-\alpha \gamma)})\geq M.
\end{align*}
Note that $w$ is radially symmetric and decreasing in the radial variable, then
\begin{align*}
w(x,t)\geq \frac{c_{3}}{c_{4}}Mt^{-\beta\gamma}  \ \mbox{for all} \ \  |x|\leq t^{\gamma}.
\end{align*}
Therefore, we obtain the conclusion.
\end{proof}

\section{\textbf{A family of functions $f$ satisfying the Osgood-type condition}}

In this section, we will construct a family of functions depending on a parameter $k>1$, which satisfy (A2) and (A3).

For $\alpha\in(1,2]$ and $k>1$, choose $\phi_{0}>\alpha^{1/(k-1)}$ and define the sequence $\phi_{i}$ by $\phi_{i+1}=\phi_{i}^{k}$.

Define $f:[0,\infty)\rightarrow [0,\infty)$ by
\begin{equation}\label{define}
f(s)= \ \left\{\begin{aligned}
&(1-\phi_{0}^{1-k})s^{k}, \ s\in J_{0}:=[0, \phi_{0}],\\
& \phi_{i}-\phi_{i-1}, \ s\in I_{i}:=(\phi_{i-1}, \phi_{i}/\alpha], \ i\geq 1,\\
&l_{i}(s), \ s\in J_{i}:=(\phi_{i}/\alpha, \phi_{i}], \ i\geq 1,
\end{aligned}\right.
\end{equation}
where $l_{i}$ denotes the linear interpolated function between the values of $f$ at $\phi_{i}/\alpha$ and $\phi_{i}$. It is clear that
$f$ satisfies (A2). For every $i\geq 1$, it is easy to see
\begin{align}\label{bound}
1<\phi_{i-1}<\phi_{i}/\alpha, \ i\geq 1,
\end{align}
and
\begin{align}\label{condition}
\lim_{i\rightarrow \infty}\phi_{i}=\infty.
\end{align}
We have
\begin{align}\label{series}
\int_{1}^{\infty}\frac{ds}{f(s)}&\geq \sum_{i}\int_{I_{i}}\frac{ds}{f(s)}\nonumber\\
&=\sum_{i=1}^{\infty}\frac{\phi_{i}/\alpha-\phi_{i-1}}{\phi_{i}-\phi_{i-1}}\nonumber\\
&=\frac{1}{\alpha}\sum_{i=1}^{\infty}\big(1-\frac{(\alpha-1)\phi_{i-1}}{\phi_{i}-\phi_{i-1}}\big)\nonumber\\
&=\frac{1}{\alpha}\sum_{i=1}^{\infty}\big(1-\frac{(\alpha-1)}{\phi_{i-1}^{k-1}-1}\big).
\end{align}
Since $\lim_{i\rightarrow \infty}\big(1-\frac{(\alpha-1)}{\phi_{i-1}^{k-1}-1}\big)=1$, by (\ref{series}), we obtain
\begin{align*}
\int_{1}^{\infty}\frac{ds}{f(s)}=\infty.
\end{align*}
Note that $f(s)$ is bounded above by $\alpha^{k}s^{k}$ for $s\geq 0$. In fact, for $s\in [0,\phi_{0}]$, it is obvious; for $s\in I_{i}=(\phi_{i-1},\phi_{i}/\alpha]$, $f(s)=\phi_{i}-\phi_{i-1}\leq \phi_{i-1}^{k}\leq s^{k}$; for $s\in J_{i}=(\phi_{i}/\alpha, \phi_{i})$, $f(s)\leq \phi_{i+1}-\phi_{i}\leq \phi_{i}^{k}\leq \alpha^{k}s^{k}$.

Define the function $\tilde{f}: [0,\infty)\rightarrow [0,\infty)$ as
\begin{equation}\label{comparison}
\tilde{f}(s)= \ \left\{\begin{aligned}
&0, \ s\in J_{0}:=[0, \phi_{0}],\\
& \phi_{i}-\phi_{i-1}, \ s\in I_{i}\cup J_{i}:=(\phi_{i-1}, \phi_{i}], \ i\geq 1,\\
\end{aligned}\right.
\end{equation}
we see that $\tilde{f}=f$ on $I_{i}$, $f\geq \tilde{f}$ on $J_{i}$. Thus $f\geq \tilde{f}$ on $[0,\infty)$.

\section{\textbf{Non-existence of local solutions}}

\begin{lemma}
If $\alpha\in(1,2]$ and $k>1+\alpha/n$, then there exists a non-negative $u_{0}\in L^{1}(\mathbb{R}^{n})$ such that (\ref{f}) has no local integral solution which is bounded in $L^{1}(\mathbb{R}^{n})$.
\end{lemma}
\begin{proof}
Since $k>1+\alpha/n$, we can choose $\beta\in (0,n)$ such that $k>(n+\alpha)/\beta$. Let $u$ be a local integral solution of (\ref{f}) with $f$ constructed as (\ref{define}) and $u_{0}=|x|^{-\beta}\chi_{R}$. Since $f\geq 0$, $u_{0}\geq 0$, $f$ is non-decreasing, from (\ref{solution}), it follows that $u(t)\geq S_{\alpha}(t)u_{0}$ for all $t\geq 0$. Therefore,
\begin{align}\label{inequality}
u(t)\geq S_{\alpha}(t)u_{0}+\int_{0}^{t}S_{\alpha}(t-s)f(S_{\alpha}(s)u_{0})ds.
\end{align}
Since $\int_{\mathbb{R}^{n}}p(t,x,y)dx=1$, we see that  $S_{\alpha}(t)$ is $L^{1}$ norm-preserving. By (\ref{inequality}) and Fubini's Theorem,
\begin{align}\label{integrable}
\|u(t)\|_{L^{1}}&=\int_{\mathbb{R}^{n}}u(t)dx\nonumber\\
&\geq \int_{\mathbb{R}^{n}}S_{\alpha}(t)u_{0}dx+\int_{\mathbb{R}^{n}}\int_{0}^{t}S_{\alpha}(t-s)f(S_{\alpha}(s)u_{0})dsdx\nonumber\\
&=\int_{\mathbb{R}^{n}}u_{0}dx+\int_{0}^{t}\int_{\mathbb{R}^{n}}f(S_{\alpha}(s)u_{0})dsdx\nonumber\\
&\geq \int_{0}^{t}\int_{\mathbb{R}^{n}}f(S_{\alpha}(s)u_{0})dxds.
\end{align}

Since $k>(n+\alpha)/\beta$, we can choose $\gamma\in (0,1/\alpha)$ such that $k>(n\gamma+1)/\beta\gamma$. Set $w(x,t)=(S(t)u_{0})(x)$. For sufficient large $i$, by (\ref{condition}), we have $\phi_{i}\geq M$ and $(c_{4}\phi_{i}/c_{3}M)^{-1/\beta\gamma}\leq t$, where $M$ is as in (\ref{m}) and $c_{3}$,  $c_{4}$ are the same constants as in (\ref{alpha}). Then by Proposition \ref{volume}, for $|x|\leq s^{\gamma}$, $0<s\leq (c_{4}\phi_{i}/c_{3}M)^{-1/\beta\gamma}$, we have $w(x,s)\geq \phi_{i}$. Since $f\geq \tilde{f}$ on $[0,\infty)$, $\phi_{i+1}=\phi_{i}^{k}$, we obtain
\begin{align}\label{integral solution}
\int_{0}^{t}\int_{\mathbb{R}^{n}}f(w(x,s))dxds&\geq \int_{0}^{(c_{4}\phi_{i}/c_{3}M)^{-1/\beta\gamma}}\int_{\mathbb{R}^{n}}f(w(x,s))dxds\nonumber\\
&\geq \int_{0}^{(c_{4}\phi_{i}/c_{3}M)^{-1/\beta\gamma}}\int_{\{x:\ |x|\leq s^{\gamma}\}}\tilde{f}(w(x,s))dxds\nonumber\\
&=\int_{0}^{(c_{4}\phi_{i}/c_{3}M)^{-1/\beta\gamma}}\int_{\{x:\ |x|\leq s^{\gamma}\}}(\phi_{i+1}-\phi_{i})dxds\nonumber\\
&\geq \frac{(\alpha-1)\phi_{i}^{k}}{\alpha}\int_{0}^{(c_{4}\phi_{i}/c_{3}M)^{-1/\beta\gamma}}\int_{\{x:\ |x|\leq s^{\gamma}\}}dxds\nonumber\\
&=\frac{(\alpha-1)\phi_{i}^{k}}{\alpha}\int_{0}^{(c_{4}\phi_{i}/c_{3}M)^{-1/\beta\gamma}}\omega_{n}s^{\gamma n}ds\nonumber\\
&=\frac{(\alpha-1)\omega_{n}}{\alpha(\gamma n+1)}\big(\frac{c_{4}}{c_{3M}}\big)^{-(rn+1)/\beta\gamma}\phi_{i}^{k-(rn+1)/\beta\gamma}.
\end{align}
By  (\ref{condition}) and note that $k>(n\gamma+1)/\beta\gamma$, we have $\phi_{i}^{k-(rn+1)/\beta\gamma}\rightarrow \infty$ as $i\rightarrow \infty$. Therefore, for such $u_{0}$, there is no integral solution of (\ref{f}) which is bounded in $L^{1}(\mathbb{R}^{n})$.
\end{proof}

\begin{theorem}
Let $\alpha\in(1,2]$ and $q\in [1,\infty)$. If $k>q(1+\alpha/n)$, then there exists a non-negative $u_{0}\in L^{q}(\mathbb{R}^{n})$ such that (\ref{f}) has no local integral solution which is in $L^{1}_{loc}(\mathbb{R}^{n})$.
\end{theorem}

\begin{proof}
Since $k>q(1+\alpha/n)$, we can choose $\beta\in (0, n/q)$ such that $k>(n+\alpha)/\beta$. Then choose $\gamma\in (0,1/\alpha)$ such that $k>(n\gamma+1)/\beta\gamma$. Let $u$ be a local integral solution of (\ref{f}) with $f$ constructed as (\ref{define}). Fix $t_{0}\in (0,1)$, for $t\in(0,t_{0})$, choose $i$ sufficiently large such that $\tilde{t}:=(c_{4}\phi_{i}/c_{3}M)^{-1/\beta\gamma}\leq t$, where $M$ is as in (\ref{m}) and $c_{3}$,  $c_{4}$ are the same constants as in (\ref{alpha}). Set $w(x,t)=(S(t)u_{0})(x)$ and let $B_{\rho}(0)$ be the ball of radius $\rho>1$ centred at 0. By (\ref{inequality}), we have
\begin{align}\label{ball}
\int_{B_{\rho}(0)}u(t)dx&\geq \int_{B_{\rho}(0)}\int_{0}^{\tilde{t}}[S_{\alpha}(t-s)f(w(\cdot,s))](x)dsdx\nonumber\\
&=\int_{0}^{\tilde{t}}\int_{B_{\rho}(0)}\int_{\mathbb{R}^{n}}p(t-s,x,y)f(w(y,s))dydxds\nonumber\\
&=\int_{0}^{\tilde{t}}\int_{\mathbb{R}^{n}}\int_{B_{\rho}(0)}p(t-s,x,y)f(w(y,s))dxdyds\nonumber\\
&\geq \int_{0}^{\tilde{t}}\int_{|y|\leq s^{\gamma}}\int_{B_{\rho}(0)}p(t-s,x,y)\frac{(\alpha-1)\phi_{i+1}}{\alpha}dxdyds\nonumber\\
&=\frac{(\alpha-1)\phi_{i+1}}{\alpha}\int_{0}^{\tilde{t}}\int_{|y|\leq s^{\gamma}}\int_{B_{\rho}(0)}p(t-s,x,y)dxdyds.
\end{align}
By (\ref{alpha}), we have
\begin{align}\label{bnrho}
\int_{B_{\rho}(0)}p(t-s,x,y)dx\geq c\int_{B_{\rho}(0)}\frac{t-s}{((t-s)^{1/\alpha}+|y-x|)^{n+\alpha}}dx.
\end{align}
Since $0<s\leq\tilde{t}<1$, then $|y|\leq s^{\gamma}<1$. The right hand integral is radially and decreasing with $|y|$. For $|x| \leq \rho \ (\rho>1)$ and $(t-s)^{-1/\alpha}>1$, choosing any unit vector $\tau$, by (\ref{bnrho}),
\begin{align}\label{bn}
\int_{B_{\rho}(0)}p(t-s,x,y)dx&\geq c\int_{B_{\rho}(\tau)}\frac{t-s}{((t-s)^{1/\alpha}+|z|)^{n+\alpha}}dz\nonumber\\
&=c\int_{B_{\rho}((t-s)^{-1/\alpha}\tau)}\frac{1}{(1+|v|)^{n+\alpha}}dv\nonumber\\
&\geq c\int_{B_{\rho}(\tau)}\frac{1}{(1+|v|)^{n+\alpha}}dv\nonumber\\
&\geq \tilde{c},
\end{align}
where $\tilde{c}$ is a positive constant. \\
By (\ref{ball}), (\ref{bn}), we obtain
\begin{align}\label{conclusion}
\|u(t)\|_{L^{1}(B_{\rho}(0))}&\geq \frac{\tilde{c}(\alpha-1)\phi_{i+1}}{\alpha}\int_{0}^{t_{0}}\int_{|y|\leq s^{\gamma}}dyds\nonumber\\
&\geq \frac{\tilde{c}(\alpha-1)\omega_{n}\phi_{i+1}}{\alpha}\int_{0}^{\tilde{t}}s^{n\gamma}ds\nonumber\\
&\geq \bar{c}\phi_{i+1}\tilde{t}^{n\gamma+1}\nonumber\\
&\geq \bar{c}\phi_{i}^{k-(n\gamma+1)/\beta\gamma}\rightarrow \infty
\end{align}
as $i\rightarrow \infty$. The proof is complete.
\end{proof}
%
%


\begin{thebibliography}{99}
\bibitem{LRS} R. Laiser, J.C. Robinson, M. Sier$\dot{\mbox{z}}$ega, Non-existence of local solutions for semilinear heat equations of Osgood type, J. Differential. Equations. 255 (2013) 3020-3028.
\bibitem{CKS}Z.-Q. Chen, K. Panki, R. Song, Dirichlet heat kernel estimates for fractional Laplacian with gradient perturbation, Ann. Probab. 40 (2012) 2483-2358.
\bibitem{TG}  T. Jakubowski, G. Serafin, Stable estimates for source solution of critical fractal Burgers equation, Nonlinear. Anal. 130 (2016) 396-407.
\bibitem{KT}  K. Bogdan, T. Byczkowski, Potential theory for the $\alpha$-stable Schr$\ddot{\mbox{o}}$dinger operators on bounded Lipschitz domains, Studia. Math. 133 (1999) 53-92.
\bibitem{Pazy} A. Pazy, Semigroups of linear operators and applications to partial differential equations. Applied Mathematical Sciences, 44. Springer-Verlag, New York, 1983.
\bibitem{Fu} H. Fujita, On the blowing up of solutions of the Cauchy problem for $u_{t}=\triangle u+u^{1+\alpha}$, J. Fac. Sci. Univ. Tokyo Sect I. 13 (1966) 109-124.
\bibitem{Fujita} H. Fujita, On some existence and nonuniqueness theorems for nonlinear parabolic equations, Proc. Symp. Pure Math., Vol.18, Part I, Amer. Math. Soc., 1968, pp. 138-161.
\bibitem{Weissler}  F.B. Weissler, Existence and nonexistence of global solutions for a semilinear heat equation, Israel J. Math. 38 (1981) 29-40.
\bibitem{QS} P. Quittner, P. Souplet, Superlinear Parabolic Problems. Blow-up, Global Existence and Steady States, Birkh$\ddot{\mbox{a}}$ Adv. Texts. Basl. Lebrb$\ddot{\mbox{u}}$, Basel, 2007.
\bibitem{HB} B. Hu, Blow-up theories for semilinear parabolic equations. Lecture Notes in Mathematics, 2018. Springer, Heidelberg, 2011.
\bibitem{PV} I. Peral,  J. L. V$\acute{\mbox{a}}$zquez, On the stability or instability of the singular solution of the semilinear heat equation with exponential reaction term, Arch. Rational Mech. Anal. 129 (1995) 201-224.
\bibitem{BC} H. Brezis, T. Cazenave, A nonlinear heat equation with singular initial data, J. Anal. Math. 68 (1996) 277-304.
\bibitem{AT} S. Abe, S. Thurner, Anomalous diffusion in view of Einsteins 1905 theory of Brownian motion.
Physica A. 356 (2005) 403-407.
\bibitem{MDH} M.M. Meerschaert, D.A. Benson, H.-P. Scheffler, B. Baeumer, Stochastic solution of space-time fractional diffusion equations, Phys. Rev. E 65 (2002) 1103-1106.
\bibitem{SJ}  M. J. Saxton, K. Jacobson, Single-particle tracking: application to membrane
dynamics. Annu. Rev. Biophys. Biomol. Struct. 26 (1997) 373-399.
\bibitem{GH}  B. Guo, Z. Huo, Global well-posedness for the fractional nonlinear Schr$\ddot{\mbox{o}}$dinger equation. Comm. Partial Differential Equations 36 (2011) 247-255.
\bibitem{CGL} Y. Cho, Z. Guo, S. Lee, A Sobolev estimate for the adjoint restriction operator, Math. Ann. 362 (2015) 799-815.
\end{thebibliography}
\end{document}